\documentclass[12pt,a4paper]{article}
\usepackage{amsfonts}
\usepackage{amsmath}
\usepackage{mathrsfs}
\usepackage{graphicx}
\usepackage{amsmath,amsthm,amssymb,amscd}
\newtheorem{Thm}{Theorem}[section]
\newtheorem{Lem}[Thm]{Lemma}
\newtheorem{Prop}[Thm]{Proposition}
\newtheorem{Def}[Thm] {Definition}

\newtheorem{Ex}  [Thm]{Example}
\newtheorem{Cor}[Thm]{Corollary}
\theoremstyle{remark}
\newtheorem{Rem} [Thm]{Remark}

\theoremstyle{claim}

\newtheorem{Que}[Thm]{Question}

\begin{document}

\begin{center} {\large \bf
 Joint Birkhoff Ergodic Average and Topological Pressure}
\end{center}
\smallskip

\begin{center}
   Xueting Tian$^{*}$\\
   \bigskip

  School of Mathematical Sciences, Fudan University, Shanghai 200433,\\ People's Republic of China \\
 \bigskip

  E-mail: {\it xuetingtian@fudan.edu.cn} \\
\end{center}

\footnotetext{ Version of   September 2014.}

\footnotetext
{$^{*}$ Tian is  supported by National Natural Science Foundation of China(No. 11301088) and  Specialized
  Research Fund for the Doctoral Program of Higher Education(No.  20130071120026).
  % and the Fundamental Research
 %Funds for the Central Universities of China(Startup Project for Young Teachers of Fudan University, No. 20520131170).
 }
\footnotetext{AMS Review:    37B10; 37B20; 37B40; 37D20; 37D25; 54H20 }

\def\abstractname{\textbf{Abstract}}

\begin{abstract}\addcontentsline{toc}{section}{\bf{English Abstract}}
  In this paper we  mainly study the dynamical complexity of Birkhoff ergodic average under the simultaneous observation of any number of continuous functions.
These results can be as   generalizations  of \cite{Barreira-Schmeling2000,To2010} etc. to study Birkhorff ergodic average  from one (or finite) observable function to any number of observable functions from  the dimensional perspective.

      For any topological dynamical  system   with  $g-$almost product property and uniform  separation property,  we show that  any {\it jointly-irregular set}(i.e., the intersection of a series of   $\phi-$irregular  sets over several  continuous functions)  either is empty or  carries  full topological pressure. In particular, if further the system is not uniquely ergodic, then the {\it completely-irregular set}(i.e., intersection of all possible {\it nonempty $\phi-$irregular} sets)  is nonempty(even forms a dense $G_\delta$ set) and carries  full topological pressure. Moreover,   {\it irregular-mix-regular sets} (i.e., intersection of some   $ \phi-$irregular sets and $ \varphi-$regular sets)   are discussed.
  Similarly, the above results are suitable for the case of   BS-dimension.

  As consequences, these results are suitable for any system such as shifts of finite type  or uniformly hyperbolic diffeomorphisms, time-1 map of uniformly hyperbolic flows, repellers, $\beta-$shifts etc..

\end{abstract}
\newpage

\section{Introduction} \setlength{\parindent}{2em}

Let $f$ be a     continuous map of a compact metric space $X$ and let $\phi:X\rightarrow \mathbb{R}$ be a continuous observable     function. A point  $x\in X$ is called to be {\it $\phi$-regular,} if the limit $$\lim_{n\rightarrow \infty}\frac1n\sum_{i=0}^{n-1}\phi(f^i(x)) $$ exists. Otherwise, $x$ is called {\it  $\phi$-irregular}(or, we say the orbit of $x$ has {\it historic behavior}). Define the $\phi$-regular set   to be the set of all $\phi$-regular points, that is, $$R(\phi,f):=\{x\in X| x\,\,is \,\,\phi-regular\};$$    and define the $\phi$-irregular set   to be the set of all $\phi$-irregular points, that is,  $$I(\phi,f):=\{x\in X| x\,\,is \,\,\phi-irregular\}.$$ By  Birkhoff's ergodic theorem,
  the $\phi-$irregular set is always of zero measure for any invariant measure.   That is, $\mu(I(\phi,f))=0$ for all invariant measure $\mu.$

Let $$I(f):=\bigcup_{\phi\in C^0(X)} I(\phi,f),$$ where
  $C^0(X)$ denotes the space of all continuous functions on $X$. $I(f)$ is called {\it irregular set}, composed of all $\phi$-irregular points for all continuous functions. Its complementary set $X\setminus I(f)$ is called {\it regular set}, denoted by $R(f).$ Every point in $R(f)$ is $\phi-$regular for all continuous functions. That is,  $$R(f)=\bigcap_{\phi\in C^0(X)}R(\phi,f).$$

% The terminology of irregular point was introduced by Ruelle in \cite{Ruelle}(called historic behavior there). As Takens said in \cite{Takens}, {\it  if this limit $\lim_{n\rightarrow \infty}\frac1n\sum_{i=0}^{n-1}\phi(f^i(x)) $ does not
%exist, it follows that `partial average' $\frac1n\sum_{i=0}^{n-1}\phi(f^i(x)) $ keep changing considerably so that
%their values give information about the epoch to which $n$ belongs: they have a history. It is
%generally believed that this is an anomaly: e.g. if the dynamics of the weather is autonomous(i.e. no climate change), then something like the average temperature should be well defined.
%There is however no justification for this belief.  }

%????add?????-------------------------------

Pesin and Pitskel \cite{Pesin-Pitskel1984} are the first to notice the phenomenon of the irregular set carrying full topological entropy in the case of the full shift on two symbols from the dimensional perspective. Barreira, Schmeling, etc. studied   the irregular set  in the setting
of shifts of finite type and beyond, see \cite{Barreira-Schmeling2000,BarBook,TV,Chen-K-Lin,Tho2012,To2010,PolWeiss,Todd,Olsen-Win} etc. Ruelle uses the terminology in \cite{Ruelle} `historic behavior' to describe irregular point and in contrast to dimensional perspective, Takens asks in \cite{Takens} for which smooth dynamical systems the points with historic behavior has positive Lebesgue measure.  Moreover,  many researchers   studied irregular set from topological or geometric viewpoint that irregular set forms dense $G_\delta$ set, see  \cite{LW,LW2,APT,Hyde,Olsen,Olsen2} etc.   So many people  paid attention to the
irregular set   and the study of irregular set is an increasingly hot topic.

In current work, we still take the dimensional viewpoint.  Firstly we recall a recent dimensional result from  \cite{To2010,Tho2012} that   which is inspired from \cite{PS} by Pfister
     and Sullivan and \cite{TV} by Takens and Verbitskiy) that

\begin{Thm}\label{Known-Result}
Let $f$ be a     continuous map of a compact metric space $X$ with (almost) specification. Then for any continuous function $\phi:X\rightarrow \mathbb{R},$  the  $\phi-$irregular set $I(\phi,f)$ either is empty or carries full topological entropy(or pressure).
\end{Thm}

Theorem  \ref{Known-Result} studied the  $\phi-$irregular set under {\bf one} observable function. 
 A natural and general question is,  how about the intersection of $\phi-$irregular sets under {\bf several} observable functions(called {\it jointly-irregular set})?  Let us recall a result of  \cite{Barreira-Schmeling2000}(Theorem 2.1)   that
\begin{Thm}\label{Bar-Sme} For systems such as topological mixing subshifts of finite type (repellers and horseshoes), the jointly-irregular set $$\bigcap_{i=1}^k I(\phi_i,f)$$ under finite observable H$\ddot{o}$lder functions $\phi_1,\phi_2,\cdots,\phi_k\,(k\geq 1)$ either is empty or carries full topological entropy.
\end{Thm}

The    special case of Theorem \ref{Known-Result}  that $I(f)=\bigcup_{\phi\in C^0(X)} I(\phi,f)$ carries full topological entropy,  was firstly  proved in \cite{Chen-K-Lin} for systems with specification property. Theorem \ref{Known-Result}  is  a refined result and Theorem \ref{Bar-Sme} is much more refined.
In present paper, we are going to give a substantial generalization of Theorem \ref{Bar-Sme}   in three directions:  (1) the system is more general; (2) the observable functions are not necessarily H$\ddot{o}$lder continuous, (3) the number of observable functions can be infinite. In particular,   Theorem \ref{Known-Result} can be generalized from one observable function to any number of observable functions. More precisely, we state the question as follows:
\begin{Que}\label{Que-Irregular-full-entropy}
Let $f$ be a     continuous map of a compact metric space $X$ with (almost) specification.  For any subset $D\subseteq C^0(X),$ whether  the  jointly-irregular set  $$\bigcap_{\phi\in D}I(\phi,f)$$  either is empty or has full topological entropy(or pressure)?
\end{Que}

We conjecture that it is true.   
To answer Question \ref{Que-Irregular-full-entropy}, in general the fundamental idea   is to construct lots of irregular points such that  you can  use  `so
many points'  to prove full topological entropy (or pressure). However, in the new setting, we need to construct lots of irregular points which are irregular for a series of functions simultaneously. So the main problem is how to deal with infinite functions. Inspired from  \cite{PS} by Pfister
     and Sullivan, we can use  variational principle to solve the problem but another condition, called uniform separation property, is required. So we can give a positive answer  in present paper under an additional assumption, uniform separation property.  Fortunately, it can be applicable to lots of classical dynamical systems, including all  uniformly hyperbolic diffeomorphisms, time-1 map of uniformly hyperbolic flows, repellers, shifts of finite type   and $\beta-$shifts etc..

Along this paper we will state and prove several theorems that give  partial positive answers to Question \ref{Que-Irregular-full-entropy}.  The results are classified in two groups:

\vspace{.3cm}

$\bullet$ {\bf The  case of finite observable functions:}
It can be dealt with similarly as in \cite{Tho2012,To2010} and it will be proved positively in  Theorem \ref{Thm-jointly-Irregular-full-entropy}. It is can be as a generalization of Theorem \ref{Bar-Sme} for systems from symbolic dynamics to more general case and for   observable functions from H$\ddot{o}$lder continuous   to just continuous.

\vspace{.3cm}

$\bullet$ {\bf The general case and the particular ``perfect" case:} Remark that if $D$ contains a constant function or a function $\varphi$ cohomologous to a constant(i.e., there is a constant $c$ and a continuous function $h$ such that $\varphi=c+h-h\circ f$), then the jointly-irregular set $$\bigcap_{\phi\in D}I(\phi,f)=\emptyset.$$
 So it is meaningful only for considering the functions with $I(\phi,f)\neq\emptyset.$ A natural and
``perfect\rq\rq  jointly-irregular set is the   intersection  of all nonempty $\phi-$irregular sets, called {\it completely-irregular set}. This set is    `minimal'  in the sense of possible non-emptiness so that if this case is true, then so does the Question \ref{Que-Irregular-full-entropy}.   For systems with $g-$almost product property and uniform separation property, we gives a positive answer(see Theorem \ref{Thm-Irregular-full-entropy-either-or} and Theorem \ref{Thm-Irregular-full-entropy}).

Remark that the case of finite observable functions does not imply the general case, because one can not find {\it finite}  observable functions such that the completely-irregular set is  a jointly-irregular set under the observation of  such    functions, see Proposition \ref{Prop-Irregular-mix-irregular}. Moreover,  we point out that it is still unknown wether the completely-irregular set can be written as a jointly-irregular set under the observation of  countable    functions.

%\bigskip

%\section{Results}

%{\bf Ref}
%It is known that irregular set is the union of all $\phi-$irregular set over all continuous functions $\phi$ and in \cite{To2010,Tho2012} Thompson showed that for systems with (almost) specification, every nonempty $\phi-$irregular set has full topological entropy.

%In this paper we give a positive answer for systems with $g-$almost product property and uniform  separation property.  That is,

\subsection{Jointly-irregular set of any number of observable functions}

 Question \ref{Que-Irregular-full-entropy} is true for systems with $g-$almost product property and uniform  separation property. Let $\mathcal{M}(X)$ denote the space of all probability measures.  Let $\mathcal{M}_f(X)$ and $\mathcal{M}^e_f(X)$   denote the space of  all $f-$invariant probability  measures and the space of all $f-$ergodic measures  respectively.

\begin{Thm}\label{Thm-Irregular-full-entropy-either-or}{\bf(  Partial Positive Answer to Question \ref{Que-Irregular-full-entropy} )}\\
Let $f$ be a     continuous map of a compact metric space $X$ with $g-$almost product property and uniform  separation property. Then for any subset $D\subseteq C^0(X),$ the the jointly-irregular  set   $\bigcap_{\phi\in D}I(\phi,f)$ either is empty or  carries  full topological pressure, that is,  $$P(\bigcap_{\phi\in D}I(\phi,f),\varphi,f)=P(X,\varphi,f)=\sup\{ {h_\mu(f)}+{\int \varphi d\mu}|\,\mu\in \mathcal{M}_f(X)\}$$ for any $\varphi\in C^0(X)$.
\end{Thm}

\begin{Rem}\label{Remark-Que-Irregular-full-entropy}
  %If    $f$   has zero topological entropy(for example, $X$ is a finite or countable set), it is trivial that the
  %jointly-irregular set  $\bigcap_{\phi\in D}I(\phi,f)$(as a subset of $X$) either is empty or has zero topological entropy.
  If $f$ is   uniquely ergodic, there is an ergodic measure $\mu$ such that for any $x\in X,$ $ \frac1{n}\sum_{i=0}^{n-1}\delta_{f^i(x)} $ converges to $\mu$ in weak$^*$ topology,   where $\delta_y$ denotes the Dirac measure supported on the point $y\in X.$  Then for any $x\in X$ and any $\phi\in C^0(X),$ $$\lim_{n\rightarrow \infty}\frac1{n}\sum_{i=0}^{n-1}\phi({f^i(x)})=\int \phi d\mu.  $$ That is, for any $\phi\in C^0(X),$ $I(\phi,f)=\emptyset.$ Thus in this case the  jointly-irregular set $\bigcap_{\phi\in D}I(\phi,f)$ is always empty.
\end{Rem}

\smallskip

 Notice that only the functions with $I(\phi,f)\neq\emptyset$ is meaningful. So we only need to consider the case of $D$ which is composed of several functions with $I(\phi,f)\neq\emptyset$. For convenience, 
 \begin{Def}\label{Def-truly-observable}({\bf Truly-observable}) 
 We call   a continuous function $\phi\in C^0(X)$ to be {\it truly-observable}, if $$I(\phi,f)\neq\emptyset.$$ That is, there is at least one orbit of some $x$ has historic behavior under the observation of $\phi.$ Otherwise, $\phi$ is called {\it trivially-observable}.
 \end{Def}    For example, any constant function is naturally trivially-observable. Roughly speaking, every truly-observable function can open an `eye' to observe irregular points.
 Let $$\hat{C}^0_f(X):=\{\phi\in C^0(X)|\,\,I(\phi,f)\neq\emptyset\}.$$ $\hat{C}^0_f(X)$ is the set of all truly-observable functions.  Remark that $$\hat{C}^0_f(X)=\emptyset\,\,\,\Leftrightarrow \,\,\,\,R(f)=X.$$ For example, the south-north map and uniquely ergodic systems are all such kind.
Let us introduce an important kind of jointly-irregular set, called completely-irregular set(or called essentially-irregular). More precisely, 
\begin{Def}\label{Def-completely-irregular}({\bf Completely-irregular set})
  define {\it completely-irregular set}  as $$CI(f):=\bigcap_{\phi\in \hat{C}^0_f(X)}I(\phi,f).$$
   In other words, the completely-irregular set $CI(f)$ is the   intersection  of all  possible nonempty $\phi-$irregular sets.  
   \end{Def}  
   Remark that it is the  `minimal'   jointly-irregular set in the sense of nonempty possibility, since if one add one more function the jointly-irregular set must be empty.  The element of completely-irregular set is called {\it completely-irregular point.}  That is, a point $x$ is completely-irregular, if $x$ is $\phi-$irregular for any $\phi\in \hat{C}^0_f(X).$

   Roughly speaking, every    completely-irregular point is persistent under all truly-observable functions. In other words, every completely-irregular point can be observed  by the `eyes' of all truly-observable functions.  Thus, if we think the regular set $R(f)$ as the `best' dynamical set, then $CI(f)$ is the `worst' dynamical set  in opposition. Remark that for any $\phi\in \hat{C}^0_f(X),$  $$CI(f)\subseteq I(\phi,f)\subseteq I(f).$$

Inspired by the  analysis of Remark \ref{Remark-Que-Irregular-full-entropy}, we only need to
 consider the system which is not uniquely ergodic.
% and has positive entropy.
Let $D\subseteq C^0(X)$ be a subset of continuous functions.
Remark that if $\bigcap_{\phi\in D}I(\phi,f)\neq \emptyset,$ then $$\bigcap_{\phi\in D}I(\phi,f)\supseteq CI(f).$$  Then Theorem \ref{Thm-Irregular-full-entropy-either-or} can be deduced from following theorem.

\begin{Thm}\label{Thm-Irregular-full-entropy} {\bf( Full Pressure of Completely-irregular set)}\\
Let $f$ be a     continuous map of a compact metric space $X$ with $g-$almost product property and uniform
separation property. Assume that   $f$ is not uniquely ergodic.
% and has positive entropy.
 Then  the  completely-irregular set   $CI( f)$ is nonempty and carries  full topological  pressure, that is,  $$P(\bigcap_{\phi\in D}I(\phi,f),\varphi,f)=P(X,\varphi,f)=\sup\{ {h_\mu(f)}+{\int \varphi d\mu}|\,\mu\in \mathcal{M}_f(X)\}$$ for any $\varphi\in C^0(X)$.
\end{Thm}

For non-uniquely ergodic systems with $g-$almost product property(or almost specification) $CI( f)$ is nonempty from  Theorem \ref{Thm-completely-Irregular-Residual} below. However,   if $f$ has positive topological entropy, then  the conclusion of full topological entropy($\varphi=0$) of Theorem \ref{Thm-Irregular-full-entropy} implies the non-emptiness of $CI( f)$ . More precisely, for $\varphi=0$,   $CI( f)$ has positive full entropy so that it is not empty.

%Theorem \ref{Thm-Irregular-full-entropy} can not be deduced from Theorem \ref{Thm-Irregular-full-entropy-either-or} directly,
 %because if $D=\hat{C}^0_f(X)$, Theorem \ref{Thm-Irregular-full-entropy-either-or} just implies $CI(f)$ either is
 % empty or has full topological entropy. However, remark that Theorem \ref{Thm-Irregular-full-entropy} implies
 % $CI(f)\neq \emptyset.$

  We will prove Theorem \ref{Thm-Irregular-full-entropy-either-or} and  Theorem \ref{Thm-Irregular-full-entropy} in Section
   \ref{section-Completely-irregular}.
%(For more description, see Theorem \ref{Thm-completely-Irregular-Residual}).

\smallskip

\subsection{Jointly-irregular set of finite observable functions}

 Notice that Theorem \ref{Thm-Irregular-full-entropy} requires the condition of uniform separation.
Here we drop the uniform separation property and try to prove a theorem as a weaker version of  Question \ref{Que-Irregular-full-entropy}. That is, we consider the  {\it jointly-irregular set} which is the intersection of {\it finite}  $\phi-$irregular sets.

\begin{Thm}\label{Thm-jointly-Irregular-full-entropy}  {\bf (Partial Positive Answer to Question \ref{Que-Irregular-full-entropy})}\\
Let $f$ be a     continuous map of a compact metric space $X$ with  (almost) specification. Then for any finite functions $\phi_1,\cdots,\phi_k\in C^0(X)(k\geq 1),$  the jointly-irregular set $\cap_{j=1}^k I(\phi_j,f)$ either is empty or carries full topological  pressure, that is,  $$P(\bigcap_{i=1}^kI(\phi_i,f),\varphi,f)=P(X,\varphi,f)=\sup\{ {h_\mu(f)}+{\int \varphi d\mu}|\,\mu\in \mathcal{M}_f(X)\}$$ for any $\varphi\in C^0(X)$.
\end{Thm}

If $f$ has uniform separation property, then this theorem can be deduced immediately from Theorem \ref{Thm-Irregular-full-entropy-either-or}. But here it is enough to assume  $f$ satisfying (almost) specification. Remark that Theorem \ref{Known-Result} is a particular case of Theorem \ref{Thm-jointly-Irregular-full-entropy}. We will prove Theorem \ref{Thm-jointly-Irregular-full-entropy} in Section \ref{section-jointly-irregular-finite-functions}.

%\section{The case of topological pressure, BS-dimension and higher-dimensional multi-fractal analysis}

\subsection{BS-dimension}\label{section-BS-dimension}

%In parallel, we remark that all  results in present paper for topological entropy can be realized for the case of
%topological pressure and BS-dimension. Let us explain more precisely. Recall that Lemma \ref{lem-PS-Estimate-Entropy}  (Variational Principle) is a characterization for topological entropy and it can deduce all above results for topological entropy. So for the case of topological pressure, it is enough just to have a Variational Principle for topological pressure. This is not difficult by slight  modification  of the proof in \cite{PS} for Lemma \ref{lem-PS-Estimate-Entropy}, for example, see \cite{PeiChen} for a detailed proof.

By the definitions of topological pressure and BS-dimension(see Section \ref{definitions} below for the definitons), it is not difficult to see that for any set $Z\subseteq X,$ the BS-dimension of $Z$ is a unique foot of Bowen's equation \begin{eqnarray}\label{eqn-pressure-dimension} P(Z,-s\varphi)=0, \,\,i.e.,\,s=BS(Z,\varphi).\end{eqnarray} So all  results in present paper for topological pressure  can be  also realized for the case of BS-dimension. For example, the conclusion of Theorem \ref{Thm-Irregular-full-entropy} can be stated

$$ BS(CI(f),f)= BS(X,f)=\sup\{\frac{h_\mu(f)}{\int \varphi d\mu}|\,\mu\in \mathcal{M}_f(X)\}$$  for any strictly positive function $\varphi\in C^0(X)$.

\subsection{Applications}

The above consequence can be applicable to lots of dynamical systems(that is, Question \ref{Que-Irregular-full-entropy} is true for these systems). For example,
    \begin{Thm}\label{Thm-Application-Irregular-full-entropy}
  For any one of following systems,  the  completely-irregular set    $CI( f)$ is not empty and carries  full topological pressure(in particular, topological entropy) and BS-dimension: \\
  {\bf(A)}.   $\,f:X\rightarrow X$ is the two-sided shift on $X=\prod^{+\infty}_{n=-\infty}Y$ where $Y=\{0,1,\cdots,k-1\}$ and   $k$ is a positive integer $\geq 2$.\\
  {\bf (A')}.   $\,f:X\rightarrow X$ is the one-sided shift on $X=\prod^{+\infty}_{n=0}Y$ where $Y=\{0,1,\cdots,k-1\}$ and   $k$ is a positive integer $\geq 2$.\\
  {\bf (B)}. $\,f:X\rightarrow X$  is a subsystem of  an Axiom A system $\,f:M\rightarrow M$ over a compact Riemannian manifold $M$ where $X$ is a hyperbolic elementary set. \\
  {\bf (B')}.  $\,f:X\rightarrow X$ is a transitive Anosov diffeomorphism of  a compact Riemannian manifold $X.$ \\
  {\bf (C)}. $f:X\rightarrow X$ is the  time-$t$ map($t\neq 0$) of a transitive Anosov flow of  a compact Riemannian manifold $X$(in this case, $f$ is partially hyperbolic).\\
  {\bf (D)}. $f:X\rightarrow X$ is   a subsystem of  a $C^1$ map  $\,f:M\rightarrow M$ over a compact Riemannian manifold $M$ where $X$ is a topological mixing and expanding invariant set(called repeller).

\end{Thm}

Remark that this result for the case of shifts of finite type can be as a generalization of \cite{Barreira-Schmeling2000}(Theorem 2.1) in three directions: (1) the observable function are not necessarily H$\ddot{o}$lder continuous, (2) the number of observable functions can be infinite and (3) topological entropy are replaced by more general concept, topological pressure. Recall a result that for $C^{1+\delta}$ conformal repellers,
it was proved in \cite{Feng} that the jointly-irregular set of finite observable functions is either empty or carries full Hausdorff dimension. In this case if $\varphi=\log \|df\|$, then $$BS(Z,\varphi)=dim_H(Z)$$ for every $Z\subseteq X.$ So     Theorem \ref{Thm-Application-Irregular-full-entropy} (D) implies that

\begin{Cor}\label{Cor-repeller-dimension}
Let $f:X\rightarrow X$ be   a subsystem of  a conformal  $C^{1+\delta}$ map  $\,f:M\rightarrow M$ over a compact Riemannian manifold $M$ where $X$ is a topological mixing and expanding invariant set(called conformal  repeller). Then for any subset $D\subseteq C^0(X),$ the the jointly-irregular  set   $\bigcap_{\phi\in D}I(\phi,f)$ either is empty or  carries  full Hausdorff dimension, that is,  $$dim_H(\bigcap_{\phi\in D}I(\phi,f))=dim_H(X).$$ 
\end{Cor}

In other words, this corollary generalizes the result of \cite{Feng} from observation of finite functions to any number of  functions.

It is known that  any system $f$ in Theorem \ref{Thm-Application-Irregular-full-entropy} is not uniquely ergodic,
has positive entropy and  satisfies specification(topological mixing $+$ shadowing property $\Rightarrow$ specification). From \cite{PS} we know that specification implies $g-$almost
product property(Proposition 2.1 in \cite{PS}). Note that the systems of (A)-(B') and (D) are all expansive. The
system $f$ of (C) is partially hyperbolic with one dimensional central bundle  and thus  $f$ is far from tangency
so that  $f$ is entropy-expansive from \cite{LiaoVianaYang}(or see \cite{DFPV,PacVie}). Recall that
from \cite{Misiurewicz} entropy-expansive implies  asymptotically $h-$expansive and from  \cite{PS}(Theorem 3.1)
any expansive or asymptotically $h-$expansive system satisfies uniform separation property.
 %It is known that expanding map is expansive and satisfies
  Thus,
  Theorem \ref{Thm-Application-Irregular-full-entropy} can be deduced from Theorem \ref{Thm-Irregular-full-entropy}.

\bigskip

Let us recall the definition of $\beta-$shift($\beta>1$) $( \Sigma_\beta, \sigma_\beta)$ in \cite{Walter}(Chapter 7.3). We only need to consider that $\beta$ is not an integer.  Consider the expansion of 1 in powers of $\beta^{-1},$ i.e. $1=\sum_{n=1}^\infty a_n\beta^{-n}$ where $a_1=[\beta]$ and $a_n=[\beta^n-\sum_{i=1}^{n-1}a_i\beta^{n-i}].$ Here $[t]$ denotes the integral part of $t\in \mathbb{R}.$ Let $k=[\beta]+1.$ Then $0\leq a_n\leq k-1$ for all $n$ so we can consider $a=\{a_n\}_1^\infty$ as a point in the space $X=\prod^{+\infty}_{n=1}Y$ where $Y=\{0,1,\cdots,k-1\}$. Consider the lexicographical ordering on $X,$ i.e. $x=\{x_n\}_1^\infty<y=\{y_n\}_1^\infty$ if $x_j<y_j$ for the smallest $j$ with $x_j\neq y_j.$ Let $f:X\rightarrow X$ denote the one-sided shift transformation. Note that $f^na\leq a$ for all $n\geq 0.$  Let $$\Sigma_\beta:=\{x=\{x_n\}_1^\infty|\,x\in X\,\text{ and }\, f^n(x)\leq a\, \text{ for all }\,n\geq 0\}.$$  Then $\Sigma_\beta$ is a closed subset of $X$ and $f(\Sigma_\beta)=\Sigma_\beta.$ Let $\sigma_\beta:=f|_{\Sigma_\beta}.$ Then $( \Sigma_\beta, \sigma_\beta)$ is one-sided $\beta-$shift. One can obtain the two-sided $\beta-$shift by letting $$\hat{\Sigma}_\beta:=\{x=\{x_n\}_{-\infty}^\infty|\,x\in \prod^{+\infty}_{n=-\infty}Y\,\text{ and }\, (x_i,x_{i+1},\cdots)\in \Sigma_{\beta}\, \text{ for all }\,i\in \mathbb{Z}\}.$$ Then $\hat{\Sigma}_\beta$ is a closed subspace of $\prod^{+\infty}_{n=-\infty}Y$ invariant under the two-sided shift $$\hat{f}:\prod^{+\infty}_{n=-\infty}Y\rightarrow \prod^{+\infty}_{n=-\infty}Y.$$
The topological entropy of $\beta-$shift($\beta>1$) is $\log \beta.$
 Remark that by Variational Principle, there is an ergodic measure with positive entropy.  Note that the Dirac measure supported on the fixed point $x=\{0\}_1^\infty\in\Sigma_\beta $ has zero entropy. So every $\beta-$shift is not uniquely ergodic.

It is known that every $\beta-$shift($\beta>1$) $( \Sigma_\beta, \sigma_\beta)$ is expansive(as a subshift of finite type) and satisfies $g-$almost product property from \cite{PS}(see the Example on P.934). So   the completely-irregular set of every $\beta-$shift  has full topological entropy($=\log\beta$). Moreover, from  \cite{Tho2012}(Lemma 5.3) every subset $Z\subseteq \Sigma_\beta$ satisfies that $$Dim_H(Z)=\frac1{log\beta}h_{top}(\sigma_\beta,Z)$$ so that $Dim_H(\bigcap_{\phi\in \hat{C}^0_f(\Sigma_\beta )}I(\phi,f))=1$. That is,

 \begin{Thm}\label{Thm-Application-beta-shift}{\bf($\beta-$shift)}
  For $\beta>1$, let $( \Sigma_\beta, \sigma_\beta)$ denote the $\beta-$shift. Then  the  completely-irregular set  $$IC(\sigma_\beta)=\bigcap_{\phi\in \hat{C}^0_f(\Sigma_\beta )}I(\phi,\sigma_\beta)$$ has full topological entropy, topological pressure,  BS-dimension and full Hausdorff dimension.

\end{Thm}

%We will discuss another kind of dimension, called BS-dimension, in subsection \ref{section-BS-dimension}.

%\subsection{}

%\subsection{Higher multi-fractal analysis of Birkhoff ergodic average}

\section{Preliminary }

\subsection{Some definitions}\label{definitions}
Firstly we recall  the definition of (almost) specification, see \cite{DGS,Sig,Bow,Bowen2,PS,To2010,Tho2012}. Let $f$ be a     continuous map of a compact metric space $X$.

\begin{Def}\label{specification}  We say that the dynamical system $f$  satisfies {\it specification property}, if the following holds:  for any $\epsilon>0$ there exists an integer $M(\epsilon)$ such that for any $k
\geq 2,$ any $k$ points $x_1,\cdots,x_k$, any  integers $$a_1\leq b_1<a_2\leq b_2\cdots<a_k\leq b_k$$ with $a_{i+1}-b_i\geq M(\epsilon)\,(2\leq i\leq k),$   there exists a point $x\in X$ such that \begin{eqnarray}\label{specification-inequality}
 d(f^j(x),f^j(x_i))<\epsilon,\,\,\,\,for \,\,a_i\leq j\leq b_i,\,1\leq i\leq k.\end{eqnarray}

\end{Def}

The original definition of specification, due to Bowen, was stronger.

\begin{Def}\label{Bowen-specification} We say that the dynamical system $f$  satisfies {\it Bowen's  specification property}, if  under the assumptions of Definition \ref{specification} and   for any  integer $p\geq M(\epsilon)+b_k-a_1,$ there exists a point $x\in X$ with $T^p(x)=x$ satisfying (\ref{specification-inequality}).

\end{Def}

We recall a result that systems with specification naturally have positive entropy. Thus if the system in  Theorem \ref{Thm-Irregular-full-entropy} has specification property, positive entropy is a natural condition. In particular, if the Bowen's specification holds, then the system is not uniquely ergodic.

\begin{Prop}\label{Remark222-Que-Irregular-full-entropy}
  Let $f$ be a     continuous map of a compact metric space $X$ with specification property. Assume $card \,X>1.$ Then $f$ has positive entropy. In particular, if the specification property is Bowen's specification, then $f$ is not uniquely ergodic.
\end{Prop}
{\bf Proof.} Positive entropy is from Proposition 21.6 of \cite{DGS}. If Bowen's specification holds, from Proposition 21.3 of \cite{DGS} we know that periodic points are dense in $X$. Notice that positive entropy implies that $X$ is uncountable. So there are two  periodic points with different orbits. Then the measures supported on the two orbits are two different periodic measures(which are all invariant and ergodic). Thus $f$ is not uniquely ergodic.  \qed

\bigskip

Recall that almost specification introduced in \cite{Tho2012} is slightly different from $g-$almost product property in \cite{PS}( Almost specification is slightly weaker). And their main ideas are same: one requires only partial shadowing of the specified orbit segments, contrary to specification property. Therefore,   in present paper we treat  {\it almost specification} same as {\it $g-$almost product property} and we only introduce the definition of $g-$almost product property as follows.  People who want to know the detailed difference, see \cite{Tho2012,PS}.  A striking and typical example of   $g-$almost product property (and almost specification) is that it applies to every $\beta-$shift\cite{Tho2012,PS}. In sharp contrast, the set of $\beta$ for which
the $\beta-$shift has specification property has zero Lebesgue measure\cite{Buzzi,Schmeling}.

Let $\Lambda_n=\{0,1,2,\cdots,n-1\}.$ The cardinality of a finite set $\Lambda$ is denoted by $|\Lambda|.$ Let $x\in X$. The dynamical ball $B_n(x,\varepsilon)$ is the set $$B_n(x,\varepsilon):=\{y\in X|\,\max\{d(f^j(x),f^j(y))|\,j\in\Lambda_n\}\leq \varepsilon\}.$$

\begin{Def}\label{blowup-function} Let $g:\mathbb{N}\rightarrow \mathbb{N}$  be a given nondecreasing unbounded map with the properties $$g(n)<n\,\,\text{ and } \lim_{n\rightarrow \infty}\frac{g(n)}n=0.$$ The function $g$ is called {\it blowup function.} Let $x\in X$ and $\varepsilon>0.$ The $g-$blowup of $B_n(x,\varepsilon)$ is the closed set
$$B_n(g;x,\varepsilon):=\{y\in X|\, \exists \Lambda\subseteq\Lambda_n ,|\Lambda_n\setminus\Lambda|\leq g(n)\,\text{ and }\,\max\{d(f^j(x),f^j(y))|\,j\in\Lambda\}\leq \varepsilon\}.$$

\end{Def}

\begin{Def}\label{product-property} We say that the dynamical system $f$  satisfies {\it $g-$almost product property}  with blowup function $g$, if  there is a nonincreasing function $m:\mathbb{R}^+\rightarrow \mathbb{N},$ such that
for any $k
\geq 2,$ any $k$ points $x_1,\cdots,x_k\in X$, any positive $\varepsilon_1,\cdots,\varepsilon_k$ and any integers $n_1\geq m(\varepsilon_1),\cdots, n_k\geq m(\varepsilon_k),$ $$\bigcap_{j=1}^k f^{-M_{j-1}}B_{n_j}(g;x_j,\varepsilon_j)\neq \emptyset,$$ where $M_0:=0,M_i:=n_1+\cdots+n_i,i=1,2,\cdots,k-1.$

\end{Def}

Now let us to recall the definition of topological entropy.

\begin{Def}\label{Def-Entropy} We say that the dynamical system $f$  satisfies {\it $g-$almost product property}  with blowup function $g$, if  there is a nonincreasing function $m:\mathbb{R}^+\rightarrow \mathbb{N},$ such that
for any $k
\geq 2,$ any $k$ points $x_1,\cdots,x_k\in X$, any positive $\varepsilon_1,\cdots,\varepsilon_k$ and any integers $n_1\geq m(\varepsilon_1),\cdots, n_k\geq m(\varepsilon_k),$ $$\bigcap_{j=1}^k f^{-M_{j-1}}B_{n_j}(g;x_j,\varepsilon_j)\neq \emptyset,$$ where $M_0:=0,M_i:=n_1+\cdots+n_i,i=1,2,\cdots,k-1.$

\end{Def}

Now we recall the definition of uniform separation property.
 Let $\mathcal{M}_x(f)$ be the set of all limits of $$\big{\{}\frac1{n}\sum_{i=0}^{n-1}\delta_{f^i(x)}\big{\}}_{n\in\mathbb{N}}$$ in weak$^*$ topology.  Set $$\mathcal{E}_n(x):=\frac1{n}\sum_{i=0}^{n-1}\delta_{f^i(x)}.$$
For $\delta>0$ and  $\varepsilon>0$, two points $x$ and $y$ are
$(\delta,n,\varepsilon)-$separated if $$|\{j:d(f^jx,f^jy)>\varepsilon,\,j\in\Lambda_n\}|\geq \delta n.$$ A subset
$E$ is  $(\delta,n,\varepsilon)-$separated if any pair of different points of $E$ are  $(\delta,n,\varepsilon)-$separated.      Let $F\subseteq \mathcal{M}(X)$ be a neighborhood of $\nu\in \mathcal{M}_f(X)$.  Define $$ X_{n,F}:=\{x\in X|\, \mathcal{E}_n(x)\in F\},$$
% \begin{eqnarray*}N(F;n,\varepsilon):=\text{maximal cardinality of a } (n,\varepsilon)-\text{separated subset of } X_{n,F},\end{eqnarray*}
and define
 \begin{eqnarray*}N(F;\delta,n,\varepsilon):=\text{maximal cardinality of a } (\delta, n,\varepsilon)-\text{separated subset of } X_{n,F}.\end{eqnarray*}
 Let $\xi=\{V_i|\,i=1,2,\cdots,k\},$ be a finite partition of measurable sets of $X$. The entropy of $\nu\in \mathcal{M}(X)$ with respect to $\xi$
 is $$H(\nu,\xi):=-\sum_{V_i\in\xi}\nu(V_i)\log \nu(V_i).$$
 We write $f^{\vee n}\xi:=\vee_{k\in \Lambda }f^{-k}\xi.$ The entropy of $\nu\in \mathcal{M}_f(X)$ with respect to $\xi$ is $$h(f,\nu,\xi):=\lim_{n\rightarrow \infty}\frac 1n H(\nu, f^{\vee n}\xi),$$ and the {\it metric entropy} of $\nu$ is $$h_\nu(f):=\sup_{\xi} h(f,\nu,\xi).$$

\begin{Def}\label{uniform-separation-property} We say that the dynamical system $f$  satisfies {\it uniform separation  property}, if  following holds. For any $\eta>0,$ there exist $\delta^*>0,\epsilon^*>0$ such  that for $\mu$ ergodic and any neighborhood $F\subseteq \mathcal{M}(X)$ of $\mu$, there exists $n^*_{F,\mu,\eta},$ such that for $n\geq n^*_{F,\mu,\eta},$ $$N(F;\delta^*,n,\epsilon^*)\geq 2^{n(h_\mu(f)-\eta)}.$$

\end{Def}

Now we recall the definition of topological pressure and entropy.

%Let $E\subseteq X,$ and $\mathfrak{F}_n(E,\epsilon)$ be the collection of all finite or countable covers of $E$
%by sets of the form $B_m(x,\epsilon)$ with $m\geq n.$ We set $$C(E;t,n,\epsilon,f):=\inf\{\sum_{B_m(x,\epsilon)\in
%\mathcal{C}}2^{-tm}:\mathcal{C}\in \mathfrak{F}_n(E,\epsilon)\},$$ and $$C(E;t,\epsilon,f):=\lim_{n\rightarrow \infty}
%C(E;t,n,\epsilon,f).$$ Then $$h_{top}(E,\epsilon,f):=\inf\{t:C(E;t,\epsilon,f)=0\}=\sup\{t: C(E;t,\epsilon,f)=\infty\}$$
%and the {\it topological entropy} of $E$ is defined as $$h_{top}(E,f);=\lim_{\epsilon\rightarrow 0}h_{top}(E,\epsilon,f).$$

% Firstly we recall the definition of topological pressure.

Let $E\subseteq X,$ $\varphi\in C^0(X)$ and $\mathfrak{F}_n(E,\epsilon)$ be the collection of all finite or countable covers of $E$ by sets of the form $B_m(x,\epsilon)$ with $m\geq n.$
 We set
 $$C(E;t,\varphi, n,\epsilon,f):=\inf\{\sum_{B_m(x,\epsilon)\in \mathcal{C}}2^{-tm+\sup_{y\in B_m(x,\epsilon)}\sum_{i=0}^{m-1}\varphi(f^i(x))}:\mathcal{C}\in \mathfrak{F}_n(E,\epsilon)\},$$
 and
 $$C(E;t,\varphi, \epsilon,f):=\lim_{n\rightarrow \infty}C(E;t,\varphi,n,\epsilon,f).$$
 Then
 $$P(E,\varphi,\epsilon,f):=\inf\{t:C(E;t,\varphi,\epsilon,f)=0\}=\sup\{t: C(E;t,\varphi,\epsilon,f)=\infty\}$$ and the {\it topological pressure} of $E$ is defined as $$P(E,\varphi,f);=\lim_{\epsilon\rightarrow 0}P(E,\varphi,\epsilon,f).$$
 In particular, if $\varphi=0,$  then the {\it topological entropy} of $E$ is defined as $$h_{top}(f,E)=P(E,0,f).$$

Now we recall BS-dimension which was introduced by Barreira and Schmeling in \cite{Barreira-Schmeling2000}.
If $\varphi$ is a strictly positive continuous function, then for each $E\subseteq X$ and each number $\epsilon>0$, define
 $$N(E;t,\varphi, \epsilon,f):=\lim_{n\rightarrow \infty}N(E;t,\varphi,n,\epsilon,f),$$ where
$$N(E;t,\varphi, n,\epsilon,f):=\inf\{\sum_{B_m(x,\epsilon)\in \mathcal{C}}2^{-t\sup_{y\in B_m(x,\epsilon)}\sum_{i=0}^{m-1}\varphi(f^i(x))}:\mathcal{C}\in \mathfrak{F}_n(E,\epsilon)\}.$$
 Set
$$BS(E,\varphi,\epsilon,f):=\inf\{t:N(E;t,\varphi,\epsilon,f)=0\}=\sup\{t: N(E;t,\varphi,\epsilon,f)=\infty\}$$ and the {\it BS-dimension} of $E$ is defined as $$BS(E,\varphi,f):=\lim_{\epsilon\rightarrow 0}BS(E,\varphi,\epsilon,f).$$

Remark that if $\varphi=1$, then $BS(E,1,f)=h_{top}(f,E).$

%\section{Completely-irregular set}\label{section-Completely-irregular}

\subsection{ Variational principle and some useful lemmas}

Firstly we recall a result from \cite{PS, PeiChen}(\cite{PS} being for topological entropy and \cite{PeiChen} being for topological pressure). We say that $f:X\rightarrow X$ is {\it saturated}, if for any $
\varphi\in C^0(X)$ and   any  compact connected nonempty set $K \subseteq \mathcal{M}_f(X),$
$$P(G_K,\varphi,f)=\inf\{h_\mu(f)+\int \varphi d\mu\,|\,\mu\in K\},$$ where $G_K=\{x\in X|\,\mathcal{M}_x(f)=K\}$(called saturated set of $K$).

%In particular, if above result only holds for $K$ to be singleton, we say $f$ is single-saturated. That is, for any
 % any $ \varphi\in C^0(X)$ and   any    $\mu\in \mathcal{M}_f(X),$
%$$P(G_\mu,\varphi,f)= h_\mu(f)+\int \varphi d\mu $$ where $G_\mu=\{x\in X|\,\mathcal{M}_x(f)=\mu\}$(called saturated set of $\mu$).

 %and  $h_{top}(f, Z)$ denotes the topological entropy of a set $Z\subseteq X.$

\begin{Lem}\label{lem-PS-Estimate-Entropy}{\bf (Variational Principle)}  Let $f$ be a     continuous map of a compact metric space $X$.
If $g-$almost product property  and uniform  separation property holds, then $f$ is saturated.
 %(2) If only $g-$almost product property  holds, then $f$ is single-saturated.

\end{Lem}

%Now let us recall a general result of upper  bound for $P(G_K,\varphi,f) $ without any assumtion from \cite{PS, PeiChen}(\cite{PS} being for topological entropy and \cite{PeiChen} being for topological pressure).

%\begin{Lem}\label{Lemma-Upper-bound-pressure} Let $f$ be a     continuous map of a compact metric space $X$.\\
%(i) Let $K\subseteq \mathcal{M}_f(X)$ be a closed subset and let $$^KG:=\{x\in X: \, \{ \mathcal{E}_n(x)\}\text{ has a limit-point in K }\}.$$ Then
%$$P(^KG,\varphi,f)\leq \sup\{h_\mu(f)+\int \varphi d\mu\,|\,\mu\in K\}.$$
%(ii) If  $K \subseteq \mathcal{M}_f(X)$ is a compact connected nonempty set, then
%$$P(G_K,\varphi,f)=\inf\{h_\mu(f)+\int \varphi d\mu\,|\,\mu\in K\}.$$ In particular, for any invariant measure
%$\mu\in \mathcal{M}_f(X),$
%$$P(G_\mu,\varphi,f)\leq h_\mu(f)+\int \varphi d\mu .$$

%\end{Lem}

\begin{Lem}\label{Lemma-C-equal} Let $f$ be a     continuous map of a compact metric space $X$.
For any $\phi\in \hat{C}^0_f(X),$ one has $$\inf_{\mu\in \mathcal{M}_f(X)}\int \phi(x)d\mu<\sup_{\mu\in \mathcal{M}_f(X)}\int \phi(x)d\mu.$$ Moreover, by Ergodic Decomposition theorem, $$\inf_{\mu\in \mathcal{M}^e_f(X)}\int \phi(x)d\mu<\sup_{\mu\in \mathcal{M}^e_f(X)}\int \phi(x)d\mu.$$

\end{Lem}
\medskip

Lemma \ref{Lemma-C-equal} is a general and direct result from the definition of $\phi-$irregular set and  can be stated more precise as one direction of following lemma.

  \begin{Lem}\label{Lemma-CompletelyIr-by-point} Let $f$ be a     continuous map of a compact metric space $X$.  Let $\phi\in C^0(X)$ and $x\in X.$ Then
 $$\phi\in \hat{C}^0_f(X),\,\,x\in  I(\phi,f) \Leftrightarrow \inf_{\mu\in \mathcal{M}_x(f)}\int \phi(x)d\mu<\sup_{\mu\in \mathcal{M}_x(f)}\int \phi(x)d\mu.$$
\end{Lem}

{\bf Proof.} On one hand, fix $\phi\in \hat{C}^0_f(X)$  and $x\in  I(\phi,f)$. By definition  there are two sequences of $n_j,m_j\uparrow +\infty$ such that the following limits exist and
$$\lim_{j\rightarrow \infty}\frac1{n_j}\sum_{i=0}^{n_j-1}\phi(f^i(x))\neq \lim_{j\rightarrow \infty}\frac1{m_j}\sum_{i=0}^{m_j-1}\phi(f^i(x)). $$ By weak$^*$ topology   one can take two convergence subsequences(if necessary) of $$\big{\{}\frac1{n_j}\sum_{i=0}^{n_j-1}\delta_{f^i(x)}\big{\}}_{n\in\mathbb{N}},\,\,\,
\big{\{}\frac1{m_j}\sum_{i=0}^{m_j-1}\delta_{f^i(x)}\big{\}}_{n\in\mathbb{N}}$$ and then  the two limits of $\mu_1$ and $\mu_2$ are in $\mathcal{M}_x(f)$ and satisfy that $$\int \phi d\mu_1=\lim_{j\rightarrow \infty}\frac1{n_j}\sum_{i=0}^{n_j-1}\phi(f^i(x))\neq \lim_{j\rightarrow \infty}\frac1{m_j}\sum_{i=0}^{m_j-1}\phi(f^i(x))= \int \phi d\mu_2.$$

On the other hand, Let  $\phi\in C^0(X)$ and $x\in X$ satisfy $$\inf_{\mu\in \mathcal{M}_x(f)}\int \phi(x)d\mu<\sup_{\mu\in \mathcal{M}_x(f)}\int \phi(x)d\mu.$$ Take two measures $\mu_1,\mu_2\in \mathcal{M}_x(f)$ such that
$$\int \phi(x)d\mu_1<\int \phi(x)d\mu_2.$$
Then  we can  take two convergence subsequences  of $$\big{\{}\frac1{n_j}\sum_{i=0}^{n_j-1}\delta_{f^i(x)}\big{\}}_{n\in\mathbb{N}},\,\,\,
\big{\{}\frac1{m_j}\sum_{i=0}^{m_j-1}\delta_{f^i(x)}\big{\}}_{n\in\mathbb{N}}$$ such the limits are $\mu_1$ and $\mu_2.$ So   $$ \lim_{j\rightarrow \infty}\frac1{n_j}\sum_{i=0}^{n_j-1}\phi(f^i(x))=\int \phi d\mu_1\neq \int \phi d\mu_2 =\lim_{j\rightarrow \infty}\frac1{m_j}\sum_{i=0}^{m_j-1}\phi(f^i(x)) .$$ Hence, $x\in I(\phi,f)$ and thus $ \phi \in \hat{C}^0_f(X).$\qed

\medskip
Remark that for systems with (almost) specification, the inverse case of Lemma \ref{Lemma-C-equal} is also true from  \cite{Tho2012,To2010}.

\medskip

\begin{Lem}\label{Lem-IC-notempty} Let $f$ be a     continuous map of a compact metric space $X$ with   (almost)  specification. Let  $\phi\in C^0 (X)$. Then    $$\inf_{\mu\in \mathcal{M}_f(X)}\int \phi(x)d\mu<\sup_{\mu\in \mathcal{M}_f(X)}\int \phi(x)d\mu\Leftrightarrow I(\phi,f)\neq \emptyset.$$

\end{Lem}
{\bf Proof} For the case of `$\Rightarrow$', see the paragraph behind of Lemma 2.1 in \cite{Tho2012}, as a corollary of Lemma 2.1 and Theorem 4.1 there, P. 5397(see Lemma 1.6 of \cite{To2010} for the case of specification).

For the case of `$\Leftarrow$', it is our above Lemma \ref{Lemma-C-equal}. \qed

\bigskip

 Recall another result that $g-$almost product  property implies entropy-density(Theorem 2.1 in \cite{PS2005}).
\begin{Lem}\label{Lemma-Entropy-dense} Let $f$ be a     continuous map of a compact metric space $X$ with   $g-$almost product  property. Then
 $f$ has entropy-dense property, that is, for any $\nu\in \mathcal{M}_f(X)$, any neighborhood $G\subseteq \mathcal{M}(X)$ of $\mu$ and any $ h_* < h_\nu(T),$ there exists an ergodic
measure $\mu\in G\cap \mathcal{M}_f(X)$ such that
%$S_\mu\neq X$ and
$h_\mu(T) > h_*.$

\end{Lem}

The following lemma is from \cite{PS}(Proposition 3.3).

\begin{Lem}\label{Lemma-Upper-Continuity}  Let $f$ be a     continuous map of a compact metric space $X$ with uniform    separation property. If the ergodic measures are entropy-dense,  then
 the entropy function
$$h_{\cdot}(f):\mathcal{M}_f(X)\rightarrow \mathbb{R},\,\mu\mapsto h_\mu(f)$$ is upper continuous.

\end{Lem}

\subsection{Cardinality of  truly-observable functions}

If one only has finite truly-observable functions, then Question \ref{Que-Irregular-full-entropy} is possibly easier to answer. However, we show that the set of truly-observable functions is uncountable for any non-uniquely ergodic system with (almost) specification.

\begin{Prop}\label{Prop-truly-obser1}  Let $f$ be a     continuous map of a compact metric space $X$ with (almost) specification. If it is not uniquely ergodic, then the set of truly-observable functions, $\hat{C}^0_f(X)$, is open and dense in $C^0(X).$
\end{Prop}

Firstly we prove a general lemma.
\begin{Lem}\label{Lem-truly-obser1}  Let $f$ be a     continuous map of a compact metric space $X$. If $\hat{C}^0_f(X)\neq \emptyset$, then $\hat{C}^0_f(X)$  is open and  dense in $C^0(X).$
\end{Lem}

{\bf Proof.} Take $\phi_0\in \hat{C}^0_f(X)$ and $x\in I(\phi_0,f).$  On one hand, we show $\hat{C}^0_f(X)$  is  dense in $C^0(X)\setminus \hat{C}^0_f(X).$ Fix $\phi\in C^0(X)\setminus \hat{C}^0_f(X).$ Then $I(\phi,f)=\emptyset$  so that $R(\phi,f)=X.$ Take $\phi_n=\frac 1n \phi_0 + \phi,\,n\geq 1.$ Then $\phi_n$ converges to $\phi$ in sup norm. By construction, it is easy to check that $x\in I(\phi_n,f),\,n\geq 1.$ That is, $\phi_n\in \hat{C}^0_f(X)$.

On the other hand, we prove that $\hat{C}^0_f(X)$ is open. Fix $\phi\in \hat{C}^0_f(X)$ and   $y\in I(\phi,f).$  Then by Lemma \ref{Lemma-CompletelyIr-by-point}, there are two different invariant measures $\mu_1,\mu_2\in \mathcal{M}_y(f)$ such that $$\int \phi d\mu_1< \int \phi d\mu_2.$$ By continuity of sup norm, we can take an open neighborhood of $\phi$, denoted by $U(\phi)$, such that for any $\varphi\in U(\phi)$ $$\int \varphi d\mu_1< \int \varphi d\mu_2.$$ Notice that $\mu_1,\mu_2\in \mathcal{M}_y(f)$.  By Lemma \ref{Lemma-CompletelyIr-by-point} $x\in I(\varphi,f),\,\forall\,\,\varphi\in U(\phi).$ This implies $\varphi\in \hat{C}^0_f(X),\,\forall\,\,\varphi\in U(\phi).$\qed

\bigskip

Now we start to prove Proposition \ref{Prop-truly-obser1}.

{\bf Proof of Proposition \ref{Prop-truly-obser1}} By Lemma \ref{Lem-truly-obser1},  we only need to prove that $\hat{C}^0_f(X)\neq \emptyset$. By assumption, there are two different invariant measures $\mu_1,\mu_2$. By weak$^*$ topology, there is a continuous function $\phi$ such that $$\int \phi d\mu_1\neq \int \phi d\mu_2.$$ By Lemma \ref{Lem-IC-notempty}, $I(\phi,f)\neq \emptyset.$ \qed

\subsection{Some topological properties of completely-irregular set}
 In a Baire space, a set is {\it residual} if it contains a countable intersection of dense open sets. Some results showed that certain
irregular sets can also be large from the topological point of view. For example,
Albeverio, Pratsiovytyi and Torbin \cite{APT}, Hyde et al \cite{Hyde} and Olsen \cite{Olsen} proved that
some kinds of irregular sets associated with integer expansion are residual. Baek
and Olsen \cite{Olsen2} discussed the set of extremely non-normal points of self-similar set
from the topological point of view. Li and Wu \cite{LW2} proved that the set of divergence
points of self-similar measure with the open set condition is either residual or empty, and they also proved in \cite{LW} that
\begin{Thm}\label{Thm-Birkhorff-Residual-Result}
Let $f$ be a     continuous map of a compact metric space $X$ with   specification. Then for any continuous function $\phi:X\rightarrow \mathbb{R},$  the  $\phi-$irregular set $I(\phi,f)$ either is empty or residual in $X$.
\end{Thm}

Now we restart to study Question \ref{Que-Irregular-full-entropy} in geometric or topological perspective and  find  that completely-irregular set is very still ``large". It can be as a generalization of Theorem \ref{Thm-Birkhorff-Residual-Result}.

\begin{Thm}\label{Thm-completely-Irregular-Residual}
Let $f$ be a  continuous map of a compact metric space $X$ with  (almost) specification. Assume that $f$ is not uniquely
ergodic. Then    the completely-irregular set $CI(f)$  is residual in $X$ 
(%or at least dense in
 $\bigcup_{\tau\in \mathcal{M}_f(X)}Support(\nu)$).
\end{Thm}

{\bf Proof.} Recall that $\hat{C}^0_f(X)\neq \emptyset$ is from Lemma \ref{Lem-truly-obser1}.

 Now we start to prove the residual property.
Firstly we give a proof for any system with Bowen's specification.
 Recall from \cite{DGS} that the set of points with maximal oscillation   $G_{max}:=\{x\in X|\,\mathcal{M}_x(f)= \mathcal{M}_f(X )\}$ is residual in $X$(Proposition 21.18 and Proposition 21.14 in \cite{DGS}).  One only needs to show that $$G_{max}\subseteq \bigcap_{\phi\in \hat{C}^0_f(X)}I_\phi.$$

 More precisely, for given $\phi\in \hat{C}^0_f(X),$ by Lemma \ref{Lemma-CompletelyIr-by-point} there are two invariant measures $\mu_1,\mu_2$ such that $$\int \phi d\mu_1 \neq  \int \phi d\mu_2.$$

 For any $z\in G_{max},$ by definition there are two subsequences  of $$\big{\{}\frac1{n}\sum_{i=0}^{n-1}\delta_{f^i(z)}\big{\}}_{n\in\mathbb{N}}$$ converging to $\mu_1,\mu_2.$   Then by  weak$^*$ topology $\int \phi d\mu_1\neq \int \phi d\mu_2$ implies  that the limit $$\lim_{n\rightarrow \infty}\frac1n\sum_{i=0}^{n-1}\phi(T^i(z)) $$ does not exist. That is, $z\in I_\phi.$ This completes the proof.

For the case of specification, the system is topologically mixing so that one can adapt the proof of \cite{DGS}(Proposition 21.18 and Proposition 21.14 in \cite{DGS})  to get that $G_{max}$ is residual in $X$. Notice that $$G_{max}\subseteq \bigcap_{\phi\in \hat{C}^0_f(X)}I_\phi$$  is a basic fact without any assumption. Therefore, Theorem \ref{Thm-completely-Irregular-Residual} is true for systems with specification. For convenience of readers, we state a rough  idea to prove that $G_{max}$ is residual in $X$. From \cite{DGS} it is a general fact that $G_{max}$ is a $G_\delta$ set. Existence of a point being in $G_{max}$ is to construct a point as follows. By entropy-density of ergodic measures(Lemma \ref{Lemma-Entropy-dense}), ergodic measures are dense in the space of invariant measures. Since  $\mathcal{M}_f(X)$ is a compact metric space, one can take a countable subset $F$, composed by ergodic measures, such that $F$ is dense in $\mathcal{M}_f(X)$. For any ergodic measure  $\nu\in F$, choose a generic point $x_\nu$ which represents the "information" of $\nu$. That is, $\mathcal{E}_n(x_\nu)$ converges to $\mu$ in weak$^*$ topology.  Then   by specification and by induction  there is a point $x$ shadowing the countable orbits of generic points more and more close(replacing the roles of periodic measures or orbits in \cite{DGS} by ergodic measures or generic points).  Then every $\nu \in F$ can be as a limit point of $\mathcal{E}_n(x) $. By density of $F$,
  $F\subseteq \mathcal{M}_x(f)$ implies $\mathcal{M}_x(f)= \mathcal{M}_f(X )$ and thus   $x$ is the needed point. Density is from topologically mixing, because   the needed shadowing point  $x$  can be chosen   to shadow  beginning from any given point in $X$.

For the case of almost specification, it is not sure that $G_{max}$ is residual in $X$. But similar as the case of specification, it is possible to  adapt the proof of \cite{DGS}(Proposition 21.18 and Proposition 21.14 in \cite{DGS}) to show that $G_{max}$ is residual in $$\Delta:=\bigcup_{\tau\in \mathcal{M}_f(X)}Support(\nu).$$ The proof of existence of a point $x$ in $G_{max}$ is similar as the case of specification. Density of $G_{max}$ in $\Delta$ is from the definition of $G_{max}$, because from weak$^*$ topology the support of any invariant measure is contained in the closure of the orbit of $x$.   Here we omit the details, see \cite{DGS} for possible modification to prove.
 \qed

\bigskip

\begin{Rem}\label{Remark-imply}
Remark that from Theorem \ref{Thm-completely-Irregular-Residual}, $CI(f)\neq \emptyset$ so that if take $D=\hat{C}^0_f(X),$ Theorem \ref{Thm-Irregular-full-entropy-either-or} implies Theorem \ref{Thm-Irregular-full-entropy}.
\end{Rem}

\bigskip

Now we study a basic property of any  completely-irregular point.
\begin{Thm}\label{Thm-completely-Irregular-Basic-Property}
Let $f$ be a  continuous map of a compact metric space $X$ with  (almost) specification. Assume that $f$ is not uniquely ergodic. Then    for any $x\in CI(f),$  $$\bigcup_{\tau\in \mathcal{M}_f(X)}Support(\nu)\subseteq \omega_f(x),$$ where $\omega_f(x)$ denotes the $\omega-$limt set of $x$. In other words, $\mu(\omega_f(x))=1$ holds for any $x\in CI(f)$ and any invariant measure $\mu\in \mathcal{M}_f(X).$ In particular, if there is an invariant measure with full support(for example, the specification is Bowen's specification), then for any $x\in CI(f),$ $\omega_f(x)=X.$
\end{Thm}

{\bf Proof.}  Assume by contradiction that there is a point $x\in CI(f),$  $$\bigcup_{\tau\in \mathcal{M}_f(X)}Support(\nu)\setminus  \omega_f(x)\neq \emptyset.$$ This implies that there exists an invariant measure $\mu$ and a point $y\in Support(\nu)\setminus \omega_f(x).$ Since $\omega_f(x)$ is a closed invariant set, we can choose $\epsilon>0$ such that the open ball centered on $y$ with radius $2\epsilon$ $B_{2\epsilon}(y)$ satisfies $$B_{2\epsilon}(y)\cap \omega_f(x)=\emptyset.$$ Then we can define a continuous function $\phi:X\rightarrow[0,1]$ such that the values of $\phi$ restricted on $\omega_f(x)$ are zero and  the values of $\phi$ restricted on the  closed ball $\bar{B}_\epsilon(y):= Closure(B_\epsilon(y))$ are 1. Note that  $y\in Support(\mu) $ implies $\mu(\bar{B}_\epsilon(y))\geq \mu(B_\epsilon(y))>0.$ Then $$\int \phi d\mu \geq \int_{\bar{B}_\epsilon(y)} \phi d\mu =\mu(\bar{B}_\epsilon(y))>0.$$ If we take an invariant measure $\nu$ supported on $ \omega_f(x)$, then $\int \phi d\nu=0.$ So $\int \phi d\nu<\int \phi d\mu.$ By Lemma \ref{Lem-IC-notempty}, $I(\phi,f)\neq \emptyset.$ So $\phi\in \hat{C}^0_f(X)$ and thus $x\in I(\phi,f).$ However, $\phi|_{\omega_f(x)}\equiv0$ implies the limit  $$\lim_{n\rightarrow\infty} \frac 1n \sum_{i=0}^{n-1}\phi(f^i(x))$$ exists and equals to 0. It contradicts to $x\in I(\phi,f).$

If the specification is Bowen's specification, from \cite{DGS} (Proposition 21.12)  a dense $G_\delta$ subset of invariant measures has support $X$. So $$\bigcup_{\tau\in \mathcal{M}_f(X)}Support(\nu)=X$$ and thus $\omega_f(x)=X.$ \qed

\section{Irregular-mix-regular   and Proof of Theorem \ref{Thm-Irregular-full-entropy-either-or}, \ref{Thm-Irregular-full-entropy} and \ref{Thm-jointly-Irregular-full-entropy}}

\subsection{Irregular-mix-regular set}
Recall that for any $\phi\in \hat{C}^0_f(X),$  $CI(f)\subseteq I(\phi,f)\subseteq I(f).$ Since the topological entropy of $CI(f)$ and $I(\phi)$ are studied, a natural interest is the complementary set $$I(\phi,f) \setminus CI(f).$$ It is the gap between completely-irregular and irregular. Remark that
\begin{eqnarray}\label{gap-irregular-and-completely}
 & &I(\phi,f)\setminus CI(f)=  I(\phi,f) \setminus (\bigcap_{\psi\in \hat{C}^0_f(X)}I(\psi,f))\nonumber\\
&=& I(\phi,f) \cap (\bigcup_{\psi\in \hat{C}^0_f(X)}R(\psi,f))=\bigcup_{ \psi\in \hat{C}^0_f(X)}(I(\phi,f)\cap R(\psi,f)).
\end{eqnarray}
Inspired by this analysis we introduce a concept called  irregular-mix-regular. A set $A\subseteq X$ is {\it irregular-mix-regular}, if $A$ is formed from  sets such as $I(\phi,f)$ and $ R(\psi,f)$ throughout  the operation  of intersection.  Firstly  we   discuss $I(\phi,f)\cap R(\psi,f)$ under two observable functions.

\begin{Thm}\label{Thm-Irregular-mix-irregular-two-function}
Let $f$ be a     continuous map of a compact metric space $X$ with $g-$almost product property and uniform  separation property. Then for any two  continuous functions $\phi,\psi\in C^0(X),$ the  irregular-mix-regular   $I(\phi, f)\cap R(\psi, f)$ either is empty or  has full topological entropy. \\
 In particular, for any $\phi \in C^0(X),$  $I(\phi,f)\setminus CI(f)$ either is empty or  has full topological pressure.
\end{Thm}

If $\psi$ is a constant function or a function cohomologous to a constant, obviously   $R(\psi,f)=X$ and thus $I(\phi, f)\cap R(\psi, f)=I(\phi, f).$
Now we give a simple   example such that $$\emptyset\neq I(\phi, f)\cap R(\psi, f) \subsetneqq I(\phi, f).$$

\begin{Ex}\label{Example-Irregular-mix-irregular}
Let $f$ be a     continuous map of a compact metric space $X$ with Bowen's specification. Fix three different periodic orbits $Orb(p_1),Orb(p_2),Orb(p_3).$ Take  $\phi$ to be a continuous function such that $$\phi|_{Orb(p_1)}=\phi|_{Orb(p_2)}=0,\,\,\,\phi|_{Orb(p_3)}=1$$ and take $\psi$ to be a continuous function such that $$\psi|_{Orb(p_3)}=\psi|_{Orb(p_2)}=0,\,\,\,\psi|_{Orb(p_1)}=1.$$
Let $\mu_1,\mu_2,\mu_3$ denote the periodic measure supported on $Orb(p_1),Orb(p_2),Orb(p_3)$ respectively. Take $K=\{\tau \mu_2+(1-\tau)\mu_3|\,\tau\in[0,1]\}.$ Then from \cite{DGS} specification implies that $G_K$ is nonempty and dense in $X,$ where $G_K=\{x\in X|\,\mathcal{M}_x(f)=K\}$(Proposition 21.14 in \cite{DGS}).  It is easy to check that $$G_K\subseteq I(\phi,f)\cap R(\psi, f)$$  so that $I(\phi,f)\cap R(\psi, f)\neq \emptyset.$  Take $K'=\{\tau \mu_1+(1-\tau)\mu_3|\,\tau\in[0,1]\}.$ Similarly, one can check that $G_{K'}\subseteq I(\phi,f) \cap I(\psi,f)$ so that $I(\phi,f)\cap R(\psi, f)\subsetneqq I(\phi, f). $

\end{Ex}

In general, we show that for any system with $g-$almost product property, every truly-observable continuous function $\phi$ satisfies   $I(\phi,f)\setminus CI(f)\neq \emptyset.$

\begin{Prop}\label{Prop-Irregular-mix-irregular}
Let $f$ be a     continuous map of a compact metric space $X$ with $g-$almost product property(or almost specification). Assume that $f$ is not uniquely ergodic. Then for any finite functions $\phi_1,\cdots,\phi_k\in \hat{C}^0_f(X)(k\geq 1),$  one has  $$\cap_{j=1}^k I(\phi_j,f)\setminus CI(f)\neq \emptyset.$$

%Then for any  $\phi\in \hat{C}^0_f(X)$ one has  $I(\phi,f)\setminus CI(f)\neq \emptyset$.\\

If further uniform separation property holds, the sets $CI(f)$ and $\cap_{j=1}^k I(\phi_j,f)\setminus CI(f)$ both have full topological entropy (deduced from  Theorem  \ref{Thm-Irregular-full-entropy} and \ref{Thm-Irregular-mix-irregular}).
\end{Prop}

\begin{Rem}\label{Remark-completely-is-nottrivial} This proposition implies $CI(f)$ is not a jointly-irregular set of finite observable functions.
\end{Rem}

{\bf Proof.} By the observation of (\ref{gap-irregular-and-completely}) we only need to find some $\psi\in \hat{C}^0_f(X)$ such that $$\cap_{j=1}^k I(\phi_j,f)\cap R(\psi,f)\neq \emptyset.$$ Before that we need the first part of Proposition 2.3 of \cite{PS2005}.

\begin{Lem}\label{Lem-PS2005Nonlinearity}  Let $f$ be a     continuous map of a compact metric space $X$ with   $g-$almost product property and $\mu$ be an ergodic measure. Then for any neighborhood $G\subseteq M(X)$ of $\mu$, there exists a closed $f$-invariant subset $Y\subseteq X$  and an integer $N_G>0$ such that for any $n\geq N_G$ and $y\in Y,$ $$\mathcal{E}_n(y) \in G.$$
\end{Lem}

Recall $\mu_1,\,\mu_2$ to be    the two ergodic measures in the proof of Theorem \ref{Thm-jointly-Irregular-full-entropy} and  for any $ 1\leq j \leq k,$
$$\int \phi_j d \mu_{1}\neq \int \phi_j d \mu_{2}.$$
By entropy density, one can choose another ergodic measure $\mu_3$ close to $\frac12(\mu_1+\mu_2)$ enough in weak$^*$ topology such that    for any $ 1\leq j \leq k,$  $$\min_{l=1,2}\int \phi_j d \mu_{l}<\int \phi_j d \mu_{3}<\max_{l=1,2}\int \phi_j d \mu_{l}.$$

Take three closed neighborhoods $G_1,G_2,G_3$ of $\mu_1,\mu_2,\mu_3$ such that  $G_1,G_2,G_3$ are pairwise disjoint and for any    measures $\nu_i\in G_i(i=1,2,3),$ $$\min_{l=1,2}\int \phi_j d \nu_{l}<\int \phi_j d \nu_{3}<\max_{l=1,2}\int \phi_j d \nu_{l}.$$
 Then by Lemma \ref{Lem-PS2005Nonlinearity} one can take three closed $f$-invariant subsets $Y_1,Y_2,Y_3\subseteq X$ and a common integer $N$ such that for any $n\geq N$ and $y\in Y_i,$ $$\mathcal{E}_n(y) \in G_i.$$ Remark that $Y_1,Y_2,Y_3$ are pairwise disjoint so that there exists a continuous function $\psi:X\rightarrow \mathbb{R}$ such that $$\psi|_{Y_1\cup Y_2}=0,\,\,\,\psi|_{Y_3}=1.$$ Take $\nu_i (i=1,2,3)$ to be three ergodic measures supported on $Y_i.$ Then by Birkhorff ergodic theorem, each $\nu_i$ is a limit point of $\mathcal{E}_n(y_i)$ for some point $y_i\in Y_i$ so that we have $\nu_i\in G_i(i=1,2,3)$. Remark  that $\nu_1,\nu_2,\nu_3$ satisfy that  $$\min_{l=1,2}\int \phi_j d \nu_{l}<\int \phi_j d \nu_{3}<\max_{l=1,2}\int \phi_j d \nu_{l}.$$

 Note that $\int \psi d\nu_1=0<1=\int \psi d \nu_3$ so that by Lemma \ref{Lem-IC-notempty}, $I(\psi,f)\neq \emptyset.$ That is, $\psi\in \hat{C}^0_f(X).$

Let $K=\{ \tau \nu_1+(1-\tau)\nu_2\}$. Then if $f$ satisfies Bowen's specification, by Proposition 21.14 in \cite{DGS} there is some  $x\in X$ such that $ \mathcal{M}_x(f)= K$(for the case of $g-$almost product property, it is not difficult to prove just by little modification to the proof of Proposition 21.14 in \cite{DGS}, here we omit the details). That $\int \phi_j d \nu_{1}\neq \int \phi_j d \nu_{2}$ implies $$\inf_{\mu\in \mathcal{M}_x(f)}\int \phi_j(x)d\mu<\sup_{\mu\in \mathcal{M}_x(f)}\int \phi_j(x)d\mu.$$ So by Lemma \ref{Lemma-CompletelyIr-by-point}, $x\in \cap_{j=1}^k I(\phi_j,f).$
Notice that $$\int \psi (x)d\nu_1=\int_{K_1} \psi (x)d\nu_1=0=\int_{K_2} \psi (x)d\nu_1=\int \psi (x)d\nu_2$$ so that $$\inf_{\mu\in \mathcal{M}_x(f)}\int \psi (x)d\mu=\sup_{\mu\in \mathcal{M}_x(f)}\int \psi(x)d\mu=0.$$ Thus by Lemma \ref{Lemma-CompletelyIr-by-point}  $x\in R(\psi,f)$. We complete the proof. \qed

\medskip

Theorem \ref{Thm-Irregular-mix-irregular-two-function} is to study irregular-mix-regular set under two observable continuous functions.
Furthermore, another interesting question is to study  the irregular-mix-regular set which is the intersection of several $\phi-$irregular sets and $\psi-$regular sets. Similar as the statements of Question \ref{Que-Irregular-full-entropy}, we ask:

\begin{Que}\label{Que-Irregular-mix-irregular}
Let $f$ be a     continuous map of a compact metric space $X$ with $g-$almost product property. Then for    any subsets $D_1,\,D_2\subseteq C^0(X),$  whether the  irregular-mix-regular set  $$(\bigcap_{\phi\in D_1}  I(\phi, f))\cap(\bigcap_{\psi\in D_2} R(\psi, f))$$ either is empty or  has full topological pressure?
\end{Que}

If $D_1\cap D_2\neq \emptyset,$ obviously $$(\bigcap_{\phi\in D_1}  I(\phi_i, f))\cap(\bigcap_{\psi\in D_2} R(\psi, f))\subseteq \bigcap_{\phi\in D_1\cap D_2}  (I(\phi, f))\cap R(\phi,f))=\emptyset.$$   If $D_1=D$ and $D_2$ is composed of a constant function, then
$$(\bigcap_{\phi\in D_1}  I(\phi, f))\cap(\bigcap_{\psi\in D_2} R(\psi, f))=\bigcap_{\phi\in D}  I(\phi, f).$$ So if Question \ref{Que-Irregular-mix-irregular} is true, then so is Question \ref{Que-Irregular-full-entropy}.

For systems with  $g-$almost product property and uniform  separation property, we have a positive answer of Question \ref{Que-Irregular-mix-irregular}  which particularly implies Theorem \ref{Thm-Irregular-mix-irregular-two-function}.

\begin{Thm}\label{Thm-Irregular-mix-irregular}
Let $f$ be a     continuous map of a compact metric space $X$ with $g-$almost product property and uniform  separation property. Then for    any subsets $D_1,\,D_2\subseteq C^0(X),$    the  irregular-mix-regular set  $$(\bigcap_{\phi\in D_1}  I(\phi, f))\cap(\bigcap_{\psi\in D_2} R(\psi, f))$$ either is empty or  has full topological pressure.
\end{Thm}

\begin{Rem}\label{Rem-similar-holds-for-lots-system}
 Remark that the results of  Theorem \ref{Thm-Irregular-mix-irregular} and Proposition \ref{Prop-Irregular-mix-irregular} are valid for the systems in Theorem \ref{Thm-Application-Irregular-full-entropy} and \ref{Thm-Application-beta-shift}.

\end{Rem}

We try to state the result of Theorem \ref{Thm-Irregular-mix-irregular} for more general topological dynamical systems.

\begin{Thm}\label{Thm-SATURATED-Irregular-mix-irregular}
Let $f$ be a     continuous map of a compact metric space $X$. If $f$ is saturated and the entropy function
$h_{\cdot}(f):\mathcal{M}_f(X)\rightarrow \mathbb{R},\,\mu\mapsto h_\mu(f)$ is upper continuous, then for    any subsets $D_1,\,D_2\subseteq C^0(X),$    the  irregular-mix-regular set  $$(\bigcap_{\phi\in D_1}  I(\phi, f))\cap(\bigcap_{\psi\in D_2} R(\psi, f))$$ either is empty or  carries  full topological pressure, that is,  $$P((\bigcap_{\phi\in D_1}  I(\phi, f))\cap(\bigcap_{\psi\in D_2} R(\psi, f)),\varphi,f)=P(X,\varphi,f)$$ for any $\varphi\in C^0(X)$.
\end{Thm}

{\bf Proof.}
 Assume
$$(\bigcap_{\phi\in D_1}  I(\phi, f))\cap(\bigcap_{\psi\in D_2} R(\psi, f))\neq \emptyset.$$ Take $x_0\in(\bigcap_{\phi\in D_1}  I(\phi, f))\cap(\bigcap_{\psi\in D_2} R(\psi, f)).$ Then by definition of $R(\psi, f),$ for any $\nu_1, \nu_2 \in \mathcal{M}_{x_0}(f)$ and $\psi\in D$, $$ \int \psi d\nu_1=\int \psi d\nu_2(=\lim_{n\rightarrow \infty}\frac1n\sum_{i=0}^{n-1}\psi(f^i(x_0))).$$ By definition of  $I(\phi, f)$ and weak$^*$ topology(see Lemma \ref{Lemma-CompletelyIr-by-point}), for any $\phi\in D_1,$ there are invariant measures $\mu_{1,\phi},\mu_{2,\phi}\in  \mathcal{M}_{x_0}(f)$ such that
$$\int \phi d \mu_{1,\phi}\neq \int \phi d \mu_{2,\phi}.$$ Note that   $$\cup_{\phi\in D_1}\{\mu_{1,\phi},\mu_{2,\phi}\}\subseteq \mathcal{M}_{x_0}(f).$$
Thus for any $\psi\in D_2$,  \begin{eqnarray}\label{pre-equal} \inf_{\phi\in D_1} \int \psi d\mu_{1,\phi}=\sup_{\phi\in D_1}\int \psi d\mu_{2,\phi}.\end{eqnarray}

Fix $\epsilon>0$ and $\varphi\in C^0(X)$. By classical Variational Principle,  we can take an ergodic measure $\mu$ such that $$h_\mu(f)+\int\varphi d\mu>P(X,\varphi,f)-\epsilon.$$ Take $\theta\in(0,1)$ close to 1 such that $$ \theta (h_\mu(f)+\int\varphi d\mu)>P(X,\varphi,f)-\epsilon$$ and $(1-\theta)\|\varphi\|<\epsilon,$ where $\|\varphi\|=\max_{x\in X}|\varphi(x)|.$  For any $ \phi\in D_1,$ if take  $\nu_{l,\phi}:=\theta\mu+(1-\theta)\mu_{l,\phi}\,(l=1,\,2)$, then  $$ h_{\nu_{l,\phi}}(f)+\int \varphi d\nu_{l,\phi}\geq \theta (h_\mu(f)+\int\varphi d\mu)-(1-\theta)\|\varphi\|>P(X,\varphi,f)-2\epsilon\,(l=1,\,2).$$ Remark that for    any $ \phi\in D_1,$
$$\int \phi d \nu_{1,\phi}\neq \int \phi d \nu_{2,\phi}$$ and for any $\psi\in D_2$, by (\ref{pre-equal}) \begin{eqnarray}\label{med-equal}\inf_{\phi\in D_1}\int \psi d\nu_{1,\phi}=\sup_{\phi\in D_1}\int \psi d\nu_{2,\phi}.\end{eqnarray}

Let $$K':= \bigcup_{\phi\in D_1}\{ \nu_{1,\phi},\,\,\nu_{2,\phi}\}$$ and let $K$ be the  convex hull of $K'.$ That is, $$K=\overline{\big{\{} \sum_{i=1}^k \lambda_i\nu_i\,\,|\,\,\,\sum_{i=1}^k \lambda_i=1,\,\lambda_i\in(0,1),\,\,\nu_i\in K'\big{\}}}.$$
Since   any
$\nu\in K'$ satisfies  $h_\nu(f)+\int\varphi d\nu>P(X,\varphi,f)-2\epsilon,$ we have $$ h_\tau(f)+\int\varphi d\tau>P(X,\varphi,f)-2\epsilon$$ for any $\tau=\sum_{i=1}^k \lambda_i\nu_i$ where $\,\lambda_i\in(0,1),\,\,\nu_i\in K',\sum_{i=1}^k \lambda_i=1.$

Since   the entropy function is upper continuous and $\rho\mapsto \int \varphi d\rho$ is continuous,  then  $$ h_\nu(f)+\int \varphi d\nu \geq  P(X,\varphi,f)-2\epsilon$$ holds for any $\nu\in K.$ Similarly, by   (\ref{med-equal}) and the continuity of functions $\psi\in D_2,$ we have \begin{eqnarray}\label{totally-equal}\inf_{\nu\in K}\int \psi d\nu=\sup_{\nu\in K}\int \psi d\nu.
\end{eqnarray}
Since  $f$ is saturated, then for the above $K$  $$P(G_K,\varphi,f)=\inf\{h_\nu(f)+\int\phi d\nu \,|\,\nu\in K\}\geq P(X,\varphi,f)-2\epsilon.$$ So the left to end the proof is only  to show  $$G_K\subseteq (\bigcap_{\phi\in D_1}  I(\phi, f))\cap(\bigcap_{\psi\in D_2} R(\psi, f)).$$ By weak$^*$ topology, it is obvious from  (\ref{totally-equal}) that  $G_K\subseteq  \bigcap_{\psi\in D_2} R(\psi, f).$

Fix $x\in G_K$ and $\phi\in D_1.$ By definition of $G_K,$ $\mathcal{M}_x(f)=K\supseteq\{ \nu_{1,\phi},\,\,\nu_{2,\phi}\}.$ Then
there are two sequences of $n_j,m_j\uparrow +\infty$ such that in weak$^*$ topology
$$\lim_{j\rightarrow \infty}\frac1{n_j}\sum_{i=0}^{n_j-1}\delta_{f^i(x)}=\nu_{1,\phi},\,\,\, \lim_{j\rightarrow \infty}\frac1{m_j}\sum_{i=0}^{m_j-1}\delta_{f^i(x)}=\nu_{2,\phi}. $$  Then
$$\lim_{j\rightarrow \infty}\frac1{n_j}\sum_{i=0}^{n_j-1}\phi(f^i(x))=\int \phi d\nu_{1,\phi} \neq \int \phi d\nu_{2,\phi}= \lim_{j\rightarrow \infty}\frac1{m_j}\sum_{i=0}^{m_j-1}\phi(f^i(x)). $$ This implies $x\in I(\phi,f).$\qed

\bigskip

{\bf Proof of Theorem \ref{Thm-Irregular-mix-irregular}. }
  By assumption, the system $f$ satisfies the conditions of Lemma \ref{Lemma-Entropy-dense} and  Lemma \ref{Lemma-Upper-Continuity}  so that  the entropy function is upper continuous. Moreover, by  Lemma \ref{lem-PS-Estimate-Entropy} $f$ is saturated. So Theorem \ref{Thm-Irregular-mix-irregular}  can be deduced from Theorem \ref{Thm-SATURATED-Irregular-mix-irregular}.\qed

%\bigskip

%\begin{Rem}\label{Rem-why-two-proofs}

\subsection{ Proof of Theorem \ref{Thm-Irregular-full-entropy-either-or} and \ref{Thm-Irregular-full-entropy}}\label{section-Completely-irregular}
  If $D_1=
 D$
 %\hat{C}^0_f(X)$
  and $D_2$ is composed of a constant function, then
$$(\bigcap_{\phi\in D_1}  I(\phi, f))\cap(\bigcap_{\psi\in D_2} R(\psi, f))=\bigcap_{\phi\in D}  I(\phi, f).$$ Therefore, Theorem \ref{Thm-Irregular-mix-irregular} implies Theorem \ref{Thm-Irregular-full-entropy-either-or}. In particular, if $D=\hat{C}^0_f(X)$,  then $$\bigcap_{\phi\in D}  I(\phi, f)=CI(f).$$  In this case $CI(f)\neq \emptyset$ by Theorem \ref{Thm-completely-Irregular-Residual}(see Remark \ref{Remark-imply}). So Theorem \ref{Thm-Irregular-full-entropy-either-or} implies Theorem \ref{Thm-Irregular-full-entropy}. \qed

%Then  can be deduced from Theorem \ref{Thm-Irregular-full-entropy} directly.

 %and then combining  Theorem \ref{Thm-completely-Irregular-Residual}(see Remark \ref{Remark-imply}) one obtains

%However, the proof of Theorem \ref{Thm-Irregular-full-entropy} in Section \ref{section-IC-Proof} directly shows full (positive) topological entropy and thus  particularly   the non-emptiness of $IC(f)$ is obvious. So we prefer to  retain the proof there.

%\end{Rem}

\subsection{Jointly-irregular set of countable observable functions}

By Theorem \ref{Thm-SATURATED-Irregular-mix-irregular}, Theorem \ref{Thm-Irregular-full-entropy-either-or} can be also stated for any system which is saturated and that entropy function is continuous.
It is still unknown whether  the continuity of entropy function in  Theorem \ref{Thm-SATURATED-Irregular-mix-irregular} and Theorem \ref{Thm-Irregular-full-entropy-either-or}  can be omitted. However,  we have a following theorem similar as  Theorem \ref{Thm-Irregular-full-entropy-either-or}  under the observation of  countable functions.

\begin{Thm}\label{Thm-SATURATED-countable-Irregular-full-entropy}
Let $f$ be a     continuous map of a compact metric space $X$.
% which is not uniquely ergodic.
%and has positive entropy.
If $f$ is saturated, then for any countable continuous functions $\phi_1,\phi_2,\cdots,\phi_n,
\cdots ,$ the the jointly-irregular  set   $\bigcap_{i=1}^\infty I(\phi_i,f)$ either is empty or  carries  full topological pressure,   that is,  $$P(\bigcap_{i=1}^\infty  I(\phi_i, f),\varphi,f)=P(X,\varphi,f)$$ for any $\varphi\in C^0(X)$.

%In particular, if all $\phi_i$ are from $\hat{C}^0_f(X)$, then $\bigcap_{i=1}^\infty  I(\phi_i, f)$ is nonempty and  $$P(\bigcap_{i=1}^\infty  I(\phi_i, f),\varphi,f)=P(X,\varphi,f)$$ for any $\varphi\in C^0(X)$.

\end{Thm}

%Recall that
%$CI(f)\neq \emptyset$ by Theorem \ref{Thm-completely-Irregular-Residual}(see Remark \ref{Remark-imply}). So,  if all $\phi_i$ are from $\hat{C}^0_f(X)$, then $\bigcap_{i=1}^\infty  I(\phi_i, f)$ contains $CI(f)$ so that must be nonempty. So we only  to deal with the first statement.

{\bf Proof.} We just need to modify the chosen $K$ in the proof of Theorem  \ref{Thm-SATURATED-Irregular-mix-irregular}.
Take $D_1=\{\phi_1,\phi_2,\cdots,\}$ and $D_2=\{1\}$.  Then $$(\bigcap_{\phi\in D_1}  I(\phi, f))\cap(\bigcap_{\psi\in D_2} R(\psi, f))=\bigcap_{i=1}^\infty  I(\phi_i, f).$$ Now assume this set is nonempty.

%Similarly as the proof of Theorem  \ref{Thm-SATURATED-Irregular-mix-irregular} take   $\mu_{1,\phi_i},\,\mu_{1,\phi_i}$ such that
%$$\int \phi_i d \mu_{1,\phi_i}\neq \int \phi_i d \mu_{2,\phi_i}.$$

Fix $\epsilon>0$ and $\varphi\in C^0(X)$. Take the same ergodic measure $\mu$ such that $$h_\mu(f)+\int\varphi d\mu>P(X,\varphi,f)-\epsilon,$$ and take $\theta\in(0,1)$ close to 1 such that $$ \theta (h_\mu(f)+\int\varphi d\mu)>P(X,\varphi,f)-\epsilon$$ and $(1-\theta)\|\varphi\|<\epsilon,$ where $\|\varphi\|=\max_{x\in X}|\varphi(x)|.$

By Lemma \ref{Lemma-C-equal}, for any $\phi_i,$ one can take an invariant measure $\mu'_{i}$ such that $$\int \phi_i d\mu\neq \int \phi_i d \mu'_{i}.$$ Then take an increasing sequence of numbers  $\theta_i\in(\theta,1)\, \uparrow$  1 such that the invariant measure  $\mu_i:=\theta_i\mu+(1-\theta)\mu'_{i}$ satisfies  $$ h_{\mu_{i}}(f)+\int \varphi d\mu_i\geq   \theta (h_\mu(f)+\int\varphi d\mu)-(1-\theta)\|\varphi\|>P(X,\varphi,f)-2\epsilon.$$ Remark that every $\mu_i$ satisfies $\int \phi d\mu \neq \int \phi d\mu_i.$

Let $$K:=\{ \mu\} \cup \bigcup_{i=1}^\infty \{t\mu_i+(1-t)\mu_{i+1}\}.$$ Since $\mu_i$ converges to $\mu,$ it  is easy to check that $K$ is connected and compact and every measure $\nu\in K$ satisfies that $h_\nu(f)+\int\varphi d\nu >P(X,\varphi,f)-2\epsilon$.

Since $f$ is saturated,   then for above $K$ one has  $$P(G_K,\varphi,f)=\inf\{h_\nu(f)+\int\varphi d\nu\,|\,\nu\in K\}\geq P(X,\varphi,f)-2\epsilon.$$ So the left to end the proof is only  to show  $$G_K\subseteq \bigcap_{i=1}^\infty I(\phi_i,f).$$

Fix $x\in G_K$ and $i\geq 1.$ By definition of $G_K,$ $\mathcal{M}_x(f)=K\supseteq\{\mu, \mu_i\}.$ Then
there are two sequences of $n_j,m_j\uparrow +\infty$ such that in weak$^*$ topology
$$\lim_{j\rightarrow \infty}\frac1{n_j}\sum_{l=0}^{n_j-1}\delta_{f^l(x)}=\mu,\,\,\, \lim_{j\rightarrow \infty}\frac1{m_j}\sum_{l=0}^{m_j-1}\delta_{f^l(x)}=\mu_i. $$  Then
$$\lim_{j\rightarrow \infty}\frac1{n_j}\sum_{l=0}^{n_j-1}\phi_i(f^l(x))=\int \phi_i d\mu \neq \int \phi_i d\mu_i= \lim_{j\rightarrow \infty}\frac1{m_j}\sum_{l=0}^{m_j-1}\phi_i(f^l(x)). $$ This implies $x\in I(\phi_i,f).$  \qed

\bigskip

  %For any $ \phi\in D_1,$ if take  $\nu_{l,\phi}:=\theta\mu+(1-\theta)\mu_{l,\phi}\,(l=1,\,2)$, then  $$ h_{\nu_{l,\phi}}(f)+\int \varphi d\nu_{l,\phi}\geq \theta (h_\mu(f)+\int\varphi d\mu)-(1-\theta)\|\varphi\|>P(X,\varphi,f)-2\epsilon\,(l=1,\,2).$$ Remark that for    any $ \phi\in D_1,$
%$$\int \phi d \nu_{1,\phi}\neq \int \phi d \nu_{2,\phi}$$ and for any $\psi\in D_2$, by (\ref{pre-equal}) \begin{eqnarray}\label{med-equal}\inf_{\phi\in D_1}\int \psi d\nu_{1,\phi}=\sup_{\phi\in D_1}\int \psi d\nu_{2,\phi}.\end{eqnarray}

%We just need to modify the chosen $K$ in the proof of Theorem \ref{Thm-SATURATED-Irregular-full-entropy} as follows.

%Fix $\epsilon>0$.
 %Let $\mu$, $\theta\in(0,1)$ and $\mu'_i=\mu'_{\phi_i}$ be the same as in Theorem \ref{Thm-SATURATED-Irregular-full-entropy} such that
% $$\theta h_\mu(f)>h_{top}(f)-\epsilon$$
% and  $$\int \phi d\mu\neq \int \phi d \mu'_{i}.$$

 % Following discussion  is similar as Theorem \ref{Thm-SATURATED-Irregular-full-entropy} by using saturated and it is left for the readers. \qed

\subsection{Proof of Theorem   \ref{Thm-jointly-Irregular-full-entropy}}\label{section-jointly-irregular-finite-functions}

%Jointly-irregular set of finite observable functions

Remark that Theorem \ref{Thm-jointly-Irregular-full-entropy} is a generalization of Theorem \ref{Known-Result}.
Here we need to consider multiple observable  functions different from  \cite{To2010,Tho2012} which only considers one function.
But the idea is just to adapt the proof of \cite{To2010,Tho2012}(Remark that the result of \cite{Tho2012} is   for the case of topological entropy but by slight modification its idea is still valid for topological pressure).   So we only give a sketch of the proof and omit the details.

%{\bf Proof of Theorem \ref{Thm-jointly-Irregular-full-entropy}.}

Assume
$$\bigcap_{j=1}^k I(\phi_j, f) \neq \emptyset.$$ Take $x\in\bigcap_{i=1}^k I(\phi_i, f).$ Then by definition of  $I(\phi_i, f)$ and weak$^*$ topology(see Lemma \ref{Lemma-CompletelyIr-by-point}), there are invariant measures $\mu_{1,j},\mu_{2,j}\in  \mathcal{M}_x(f)(1\leq j\leq k)$ such that
$$\int \phi_j d \mu_{1,j}\neq \int \phi_j d \mu_{2,j},\,\,1\leq j \leq k.$$

For a fixed $1\leq j \leq k,$  obviously  the solutions $(\lambda_1,\cdots,\lambda_k)$ of linear equation  $$ \sum_{i=1}^k(\int \phi_j d \mu_{1,i}-\int \phi_j d \mu_{2,i})\lambda_i=0$$ form a $k-1$ dimensional  closed linear subspace of $\mathbb{R}^k$, denoted by $P_j.$  Notice that $$\bigcup_{1\leq j\leq k}P_j$$ is the union of finite $k-1$ dimensional  closed linear subspaces so that  $$\mathbb{R}^k\setminus (\bigcup_{1\leq j\leq k}P_j)$$ is open and dense in $\mathbb{R}^k.$ One can take $k$ positive numbers of $\lambda_1,\cdots,\lambda_k$ such that for any $ 1\leq j \leq k,$
$$\sum_{i=1}^k(\int \phi_j d \mu_{1,i}-\int \phi_j d \mu_{2,i})\lambda_i\neq 0.$$ For example, if $k=3$, $(\lambda_1,\lambda_2,\lambda_3)$ is chosen  in the first octant of $\mathbb{R}^3$.

Let $\theta_j=\frac{\lambda_j}{\sum_{i=1}^k \lambda_i}.$ Then all $\theta_j$ are positive and $\sum_{j=1}^k\theta_j=1.$ If we define $$\mu_l= \sum_{i=1}^k \theta_i\mu_{l,i},\,\,\,l=1,\,2,$$ then for    any $ 1\leq j \leq k,$
$$\int \phi_j d \mu_{1}\neq \int \phi_j d \mu_{2}.$$

Fix $\epsilon>0$ and $\varphi\in C^0(X)$. By classical Variational Principle,  we can take an ergodic measure $\mu$ such that $$h_\mu(f)+\int\varphi d\mu>P(X,\varphi,f)-\epsilon.$$ Take $\theta\in(0,1)$ close to 1 such that $$ \theta (h_\mu(f)+\int\varphi d\mu)>P(X,\varphi,f)-\epsilon$$ and $(1-\theta)\|\varphi\|<\epsilon,$ where $\|\varphi\|=\max_{x\in X}|\varphi(x)|.$ Then the two invariant measures  $\nu_l:=\theta\mu+(1-\theta)\mu_{l}\,(l=1,\,2)$ satisfy   $$ h_{\nu_{l}}(f)+\int \varphi d\nu_{l}\geq \theta (h_\mu(f)+\int\varphi d\mu)-(1-\theta)\|\varphi\|>P(X,\varphi,f)-2\epsilon\,(l=1,\,2).$$ Remark that for    any $ 1\leq j \leq k,$
$$\int \phi_j d \nu_{1}\neq \int \phi_j d \nu_{2}.$$

By \cite{PS2005}, when $f$ has the almost specification, we can find  two sequence of ergodic measures $\nu_{l,i}\in \mathcal{M}_f(X)$ such that $$h_{\nu_{l,i}}(f)\rightarrow h_{\nu_l},\,\,{\text and}\,\,\,\nu_{l,i} \rightarrow \nu_l\,(l=1,2)$$ in weak$^*$ topology. Therefore, we can take two measures belonging to these two sequence which we called $\rho_1$ and $\rho_2$ respectively such that $$ h_{\rho_{l}}(f)+ \int \varphi d\rho_l>P(X,\varphi,f)-2\epsilon\,(l=1,\,2)$$ and for    any $ 1\leq j \leq k,$
$$\int \phi_j d \rho_{1}\neq \int \phi_j d \rho_{2}.$$  This is the first step of \cite{To2010,Tho2012} to choose two good ergodic measures but crucial because  the ergodicity is  important in the proof of \cite{To2010,Tho2012} to avoid the use of uniform separation.

Then we can follow the proof of \cite{To2010,Tho2012} to complete the proof. Roughly speaking, using the above two ergodic measures to construct a set $F\subseteq \cap_{j=1}^k I(\phi_j,f)$ such that the topological pressure of $F$ is larger than $P(X,\varphi,f)-2\epsilon.$ In this process, Entropy Distribution Principle plays an important role. One can see \cite{To2010,Tho2012} for more details. \qed

\section{Simple Applications}

In this section we compute a detailed example as a simple application.

Let $f:S^1\rightarrow S^1, x\mapsto 2x \text{ mod } 1$. It is known that $f$ is expansive and satisfies Bowen's specification(it is also known $f$ is topologically conjugated to the shift of two symbols except at countably many points).   Consider   $$D:=\bigcup_{n\geq 1}(\{sin2n\pi x\}\cup \{cos2n\pi x\}).$$ $D$ is a particular subset of $C^0(S^1)$ consisted of some sine functions and
  cosine functions from the basis of Fourier series. and it is easy to check every function in $D$  is truly-observable. In fact, let $\mu$ be the Dirac measure supported on the fixed point $0\in S^1$ and let $\mu_n$ be the periodic measure supported on the periodic orbit which is the $\omega-$limit set of the orbit $$O_n:=\{ \frac{1}{7n},\,\frac{2}{7n},\,\frac{4}{7n},\,\cdots,\,\frac{2^k}{7n},\,\cdots \}\,mod\,1.$$
Then $$\int sin 2n\pi(x) d\mu=\lim_{k\rightarrow \infty}\frac1k \sum_{j=0}^{k-1} sin 2n\pi 0=0 $$ and
$$\int sin 2n\pi(x) d\mu_n=\lim_{k\rightarrow \infty}\frac1k \sum_{j=0}^{k-1} sin 2n\pi \frac{2^j}{7n}=\frac 13(sin \frac27\pi+sin \frac47\pi+sin \frac87\pi)\neq 0 .$$
Similarly, we can check that $$\int cos 2n\pi(x) d\mu=\lim_{k\rightarrow \infty}\frac1k \sum_{j=0}^{k-1} cos 2n\pi 0=1 $$ and
$$\int cos 2n\pi(x) d\mu_n=\lim_{k\rightarrow \infty}\frac1k \sum_{j=0}^{k-1} cos 2n\pi \frac{2^j}{7n}=\frac 13(cos \frac27\pi+cos \frac47\pi+cos \frac87\pi)\neq 1. $$

By Lemma \ref{Lem-IC-notempty}, $D\subseteq \hat{C}^0_f(X).$ So    $$\bigcap_{\phi\in D}I(\phi,f)$$ carries full topological entropy and full Hausdorff dimension.  In other words, the jointly-irregular set under the observation of above  infinite sine functions and
  cosine functions is still meaningful  and carries full topological pressure, topological entropy, BS-dimension and full Hausdorff dimension.

%Remark that in \cite{PS} almost specification is called $g-$almost product property.

%Remark that the results of  Theorem \ref{Thm-jointly-Irregular-full-entropy} are valid for the systems in Theorem \ref{Thm-Application-Irregular-full-entropy} and \ref{Thm-Application-beta-shift}.

\end{document}